\documentstyle[12pt]{article}
       \font\tenmsb=msbm10
       \font\sevenmsb=msbm7
       \font\fivemsb=msbm5
       \catcode`\@=11
       \ifx\amstexloaded@\relax\catcode`\@=\active
       \endinput\else\let\amstexloaded@\relax\fi

       \newfam\msbfam
       \textfont\msbfam=\tenmsb
       \scriptfont\msbfam=\sevenmsb
       \scriptscriptfont\msbfam=\fivemsb

\newtheorem{Theorem}{Theorem}[section]

\newtheorem{Corollary}{Corollary}[section]
\newtheorem{Remark}{Remark}[section]
\newtheorem{Example}{Example}[section]
\newtheorem{Definition}{Definition}[section]

\@addtoreset{equation}{section}

\newcommand{\proof}{\par\medbreak\it Proof. \rm}

\begin{document}
\setlength{\columnsep}{5pt}
\title{On the perturbation and expression for the core inverse of linear operator in Hilbert space$^{\star}$ }
\par\vspace{0.5in}
\author{\small  Qianglian  ${\rm Huang}_,^\dagger$\quad Saijie Chen,\quad  Lanping Zhu \\
 {\small  \it{ School of Mathematical Sciences, Yangzhou University,
 Yangzhou 225002, China}}}
\date{}
\maketitle \footnotetext{\scriptsize{$^{\star}$This research was supported by the National Natural Science Foundation of China (11771378) and the Yangzhou University Foundation (2016zqn03) for Young Academic Leaders.}} \footnotetext{\scriptsize{ $^{\dagger}$Corresponding author. {\it E-mail:} huangql@yzu.edu.cn; 2428245141@qq.com; lpzhu@yzu.edu.cn.}}
 \maketitle\vspace*{-1\baselineskip}
\par
 \noindent{\bf Abstract}
  {\small \
In this note, we provide some sufficient and necessary conditions for  the
core inverse of the perturbed operator to have the simplest possible expression. The results improve the recent work by H. Ma (Optimal perturbation bounds for the core inverse, Appl. Math.
Comput. 336 (2018) 176-181.).
\par\vspace{0.1in} \noindent{\bf{\it AMS classification}}: \ 15A09; \
47A05; \ 65F20
\vspace*{-0.5\baselineskip}
\section{Introduction and Preliminaries}
\indent\quad\ Let $X$  be a Hilbert space and $B(X)$ denote the Banach
space of all bounded linear operators from $X$ into itself. For any $T\in B(X)$, we denote by $N(T)$
and $R(T)$ the null space and respectively, the range of
$T$. The identity operator will be denoted by $I$.
\par Recall that an operator $S\in B(X)$ is said to be a generalized inverse of
$T\in B(X)$ if $S$ satisfies: \vspace*{-0.5\baselineskip}
$$\quad (1)\ TST=T \quad {\rm and}\quad (2)\ STS =S.$$
The generalized inverse is not unique in general. To force its uniqueness,
some further conditions must be imposed. The perhaps most convenient additional conditions are \vspace*{-0.3\baselineskip}
$$ (3)\ (TS)^*=TS;\ (4)\
(ST)^*=ST;\ (5)\  TS=ST; \ (6)\ ST^2=T \ {\rm and} \ (7)\ TS^2=S.$$
\quad\ Let
$\theta\subset \{1,2,3,4,5,6,7\}$ be a nonempty set. If $S\in B(X)$ satisfies the equation $(i)$
for all $i\in \theta$, then $S$ is said to be a $\theta$-inverse of
 $T$, which is  denoted  by $T^\theta$. As is well-known, each kind of $\theta$-inverse has
its own specific characteristics \cite{CW,RDD}. For example, the $\{1\}$-inverse is  inner
inverse,  $\{2\}$-inverse is outer
inverse and
$\{1,2\}$-inverse of $T$ is   generalized
inverse.  The
$\{1,2,3,4\}$-inverse is exactly the Moore-Penrose inverse
and $\{1,2,5\}$-inverse is  group inverse.  \par As an important generalized inverse in some sense in-between Moore-Penrose inverse and group inverse \cite{BaT1,BaT2,RDD}, the core inverse is first introduced by Baksalary and Trenkler in the matrix case.\vspace*{-0.5\baselineskip}
\begin{Definition}\label{Def1.1} ${\rm\cite{BaT1}}$ Let $T\in C^{n\times n}$ and ${\rm ind}\ T\leq 1$. A matrix $T^{\tiny\textcircled{\tiny\#}} \in C^{n\times n}$ satisfying
 \vspace*{-0.3\baselineskip} $$TT^{\tiny\textcircled{\tiny\#}}= P_T \quad{\rm and}\quad  R(T^{\tiny\textcircled{\tiny\#}} )\subseteq R(T)$$
 is called the core inverse of $T$, where ${\rm ind}\ T$  is the smallest nonnegative integer $k$ satisfying ${\rm Rank}\ T^k={\rm Rank}\ T^{k+1}$ and $P_T$ is the orthogonal projector on $R(T)$.
\end{Definition}
 \vspace*{-0.4\baselineskip}\quad\  Wang and Liu \cite{WaL}  proved that the core inverse $T^{\tiny\textcircled{\tiny\#}}$ is the unique matrix satisfying equations (1), (3) and (7). This means that the core inverse $T^{\tiny\textcircled{\tiny\#}}$ is exactly the $\{1,3,7\}$-inverse in the matrix case. If we define core inverse of an operator $T\in B(X)$ in the same way, then we have the
problem because the index of $T$ is defined by the rank of the matrix. In \cite{RDD}, Raki${\rm\acute{c}}$, Din${\rm\check{c}}$i${\rm\acute{c}}$ and Djordjevi${\rm\acute{c}}$ introduced the definition of core inverse in $B(X)$ and proved that it is equivalent to Definition 1.1 in the matrix case.
 \vspace*{-0.3\baselineskip} \begin{Definition}\label{Def1.2} ${\rm\cite{RDD}}$  An operator $T^{\tiny\textcircled{\tiny\#}}\in B(X)$ is core inverse of $T\in B(X)$ if \vspace*{-0.5\baselineskip}
 $$(1)\ TT^{\tiny\textcircled{\tiny\#}} T=T; \quad  R(T^{\tiny\textcircled{\tiny\#}} )=R(T)\quad and\quad N(T^{\tiny\textcircled{\tiny\#}} )=N(T^*).$$
 \end{Definition}
 \vspace*{-0.5\baselineskip}\quad\ In the same paper, they proved that the core inverse in infinite dimensional Hilbert space is precisely the $\{1,2,3,6,7\}$-inverse and emphasized that none of the equations (1), (2), (3), (6) and (7) can be removed.
  \vspace*{-0.3\baselineskip}\begin{Theorem}\label{The1.1}${\rm\cite{RDD}}$
 The operator $T^{\tiny\textcircled{\tiny\#}}\in B(X)$ is core inverse of $T\in B(X)$ if and only if $T^{\tiny\textcircled{\tiny\#}}$ satisfies equations (1), (2), (3), (6) and (7).
\end{Theorem}
 \vspace*{-0.3\baselineskip}\quad\ The perturbation and expression of the core inverse in matrix case is recently studied by Ma \cite{M}, which is inspired by the following well known fact: if $T\in C^{n\times n}$ is invertible and $T^{-1}\in C^{n\times n}$ is its inverse, then for any $\delta T\in C^{n\times n}$ satisfying
 $\|T^{-1}\delta T\|<1$, $\overline{T}=T+\delta T$ is invertible and its inverse is $\overline{T}^{-1}=T^{-1}(I+\delta T T^{-1})^{-1}=(I+T^{-1}\delta T)^{-1}T^{-1}$.
 \vspace*{-0.3\baselineskip} \begin{Theorem}\label{The1.2} ${\rm\cite{M}}$ Let $T\in C^{n\times n}$ be a matrix
with a  core inverse $T^{\tiny\textcircled{\tiny\#}}\in C^{n\times n}$ and
$\delta T\in C^{n\times n}$ with $\|T^{\tiny\textcircled{\tiny\#}}\delta T\|<1$. If $TT^{\tiny\textcircled{\tiny\#}}\delta T=\delta T$, then \vspace*{-0.5\baselineskip}
$$B=T^{\tiny\textcircled{\tiny\#}}(I+\delta TT^{\tiny\textcircled{\tiny\#}})^{-1}=(I+T^{\tiny\textcircled{\tiny\#}}\delta T)^{-1}T^{\tiny\textcircled{\tiny\#}}$$ \vspace*{-0.5\baselineskip}  is a  core inverse of $\overline{T}=T+\delta T$, $\overline{T}B=TT^{\tiny\textcircled{\tiny\#}}$ and  \vspace*{-0.1\baselineskip}
 $$\frac{\|T^{\tiny\textcircled{\tiny\#}}\|}{1+\|T^{\tiny\textcircled{\tiny\#}}\delta T\|}\leq \|B\|\leq\frac{\|T^{\tiny\textcircled{\tiny\#}}\|}{1-\|T^{\tiny\textcircled{\tiny\#}}\delta T\|}.$$
\end{Theorem}
 \vspace*{-0.3\baselineskip}\quad\ In the next section, we extend Theorem 1.2 to the case of bounded linear operators in Hilbert space and prove that the
core inverse of the perturbed operator has the simplest possible expression if and only if the perturbation is range-preserving. We give a direct proof since it is not only short, but also constructive. We also propose two interesting examples to illustrate our results in Section 2 and conclude in Section 3.
\section{Main Results}
 \vspace*{-0.5\baselineskip}\quad\ Since the core inverse is exactly the $\{1,3,7\}$-inverse in the matrix case \cite{WaL}, we first consider the $\{1,3,7\}$-inverse.
 \begin{Theorem}\label{The2.1} Let $X$  be a Hilbert
 space. Let $T\in B(X)$
with a $\{1,3,7\}$-inverse $T^{\{1,3,7\}}\in B(X)$ and
$I+T^{\{1,3,7\}}\delta T: X\rightarrow X$ be invertible  with $\delta T\in B(X)$. Then
$$B=T^{\{1,3,7\}}(I+\delta TT^{\{1,3,7\}})^{-1}=(I+T^{\{1,3,7\}}\delta T)^{-1}T^{\{1,3,7\}}$$  is a $\{1,3,7\}$-inverse of $\overline{T}=T+\delta T$ if and only if \vspace*{-0.5\baselineskip}
$$R(\overline{T})\subseteq R(T)$$ \vspace*{-0.5\baselineskip}
if and only if \vspace*{-0.5\baselineskip}
$$R(\overline{T})= R(T)$$\vspace*{-0.5\baselineskip}
if and only if
$$ TT^{\{1,3,7\}}\overline{T}= \overline{T} \quad or \quad TT^{\{1,3,7\}}\delta T= \delta T. $$
\end{Theorem}
\quad\ {\it In this case, $\overline{T}B=TT^{\{1,3,7\}}$, $\|B\|\leq \|T^{\{1,3,7\}}\|\cdot\|(I+T^{\{1,3,7\}}\delta T)^{-1}\|$ and }$$ \|B-T^{\{1,3,7\}}\|\leq \|T^{\{1,3,7\}}\|\cdot\|(I+T^{\{1,3,7\}}\delta T)^{-1}\|\cdot\|T^{\{1,3,7\}}\delta T\|.$$
\proof It follows from the spectral theory that $I+\delta T T^{\{1,3,7\}}: X\rightarrow X$ is also invertible. We can easily check \vspace*{-0.5\baselineskip}
$$T^{\{1,3,7\}}(I+\delta TT^{\{1,3,7\}})^{-1}=(I+T^{\{1,3,7\}}\delta T)^{-1}T^{\{1,3,7\}}$$
and so $B$ is well defined. If $B$  is a $\{1,3,7\}$-inverse of $\overline{T}$, then $(\overline{T}B)^2=\overline{T}B\overline{T}B=\overline{T}B$ and $(\overline{T}B)^*=\overline{T}B$. Hence \vspace*{-1\baselineskip}
 $$X=R(\overline{T}B)\dot{+}N(\overline{T}B)=R(\overline{T})\dot{+} N(\overline{T}B),$$
 \vspace*{-0.3\baselineskip} where $\dot{+}$ denotes the orthogonal topological direct sum. Similarly, \vspace*{-0.1\baselineskip} $$X=R(TT^{\{1,3,7\}})\dot{+}N(TT^{\{1,3,7\}})=R(T)\dot{+}N(TT^{\{1,3,7\}}).$$
\quad\quad Now we can claim that $N(TT^{\{1,3,7\}})\subseteq N(\overline{T}B).$ In fact, for any $y\in N(TT^{\{1,3,7\}})$, we get $T^{\{1,3,7\}}y\in N(T)$ and
\begin{eqnarray*}
0&=&(B\overline{T}-I)B\overline{T}T^{\{1,3,7\}}y\\
 &=&(I+T^{\{1,3,7\}}\delta T)^{-1}(T^{\{1,3,7\}}\overline{T}-I-T^{\{1,3,7\}}\delta T)B\overline{T}T^{\{1,3,7\}}y\\
 &=&(I+T^{\{1,3,7\}}\delta T)^{-1}(T^{\{1,3,7\}}T-I)B\overline{T}T^{\{1,3,7\}}y.
 \end{eqnarray*}
Then $(T^{\{1,3,7\}}T-I)B\overline{T}T^{\{1,3,7\}}y=0$ and so
\begin{eqnarray*}
0&=&(T^{\{1,3,7\}}T-I)B\overline{T}T^{\{1,3,7\}}y\\
&=&(T^{\{1,3,7\}}T-I)B\delta TT^{\{1,3,7\}}y\\
 &=&(T^{\{1,3,7\}}T-I)B(I+\delta TT^{\{1,3,7\}}-I)y\\
  &=&(T^{\{1,3,7\}}T-I)T^{\{1,3,7\}}y-(T^{\{1,3,7\}}T-I)By\\
 &=&-T^{\{1,3,7\}}y+(I-T^{\{1,3,7\}}T)By,
 \end{eqnarray*}
 i.e., $(I-T^{\{1,3,7\}}T)By=T^{\{1,3,7\}}y$.  Hence we obtain
 \begin{eqnarray*}
\overline{T}By&=&\overline{T}B\overline{T}By\\
&=&\overline{T}(I+T^{\{1,3,7\}}\delta T)^{-1}T^{\{1,3,7\}}(T+\delta T)By\\
 &=&\overline{T}(I+T^{\{1,3,7\}}\delta T)^{-1}(T^{\{1,3,7\}}T+T^{\{1,3,7\}}\delta T)By\\
  &=&\overline{T}(I+T^{\{1,3,7\}}\delta T)^{-1}[(T^{\{1,3,7\}}T-I)+(I+T^{\{1,3,7\}}\delta T)]By\\
 &=&\overline{T}(I+T^{\{1,3,7\}}\delta T)^{-1}(-T^{\{1,3,7\}}y)+\overline{T}By\\
 &=&-\overline{T}By+\overline{T}By=0,
 \end{eqnarray*}
 which  implies $y\in N(\overline{T}B)$. Therefore, $N(TT^{\{1,3,7\}})\subseteq N(\overline{T}B)$ and
$$R(\overline{T})=[N(\overline{T}B)]^\perp\subseteq [N(TT^{\{1,3,7\}})]^\perp=R(T).$$
\par Conversely, if $ R(\overline{T})\subseteq R(T)$, then for all $x\in X$, \vspace*{-0.5\baselineskip}
$$ \overline{T}x\in R(\overline{T})\subseteq R(T)=R(TT^{\{1,3,7\}})=N(I-TT^{\{1,3,7\}})$$ \vspace*{-0.5\baselineskip}
and so $ (I-TT^{\{1,3,7\}})\overline{T}x=0$, that is, $\overline{T}=TT^{\{1,3,7\}}\overline{T}$. Hence
$$\delta T=\overline{T}-T=TT^{\{1,3,7\}}\overline{T}-TT^{\{1,3,7\}}T=TT^{\{1,3,7\}}\delta T$$ \vspace*{-0.5\baselineskip} and \vspace*{-0.5\baselineskip}
$$\overline{T}=T+\delta T=T+TT^{\{1,3,7\}}\delta T=T(I+T^{\{1,3,7\}}\delta T).$$
\quad\ Noting that $I+T^{\{1,3,7\}}\delta T$ is invertible, we can obtain $R(\overline{T})= R(T)$. In the following, we shall show that
$B$ is a $\{1,3,7\}$-inverse of $\overline{T}$. In fact, \vspace*{-0.2\baselineskip}
$$\overline{T}B\overline{T}=T(I+T^{\{1,3,7\}}\delta T)(I+T^{\{1,3,7\}}\delta T)^{-1}T^{\{1,3,7\}}\overline{T}=TT^{\{1,3,7\}}\overline{T}=\overline{T},$$ \vspace*{-0.5\baselineskip}
and  \vspace*{-0.5\baselineskip}$$\overline{T}B=T(I+T^{\{1,3,7\}}\delta T)(I+T^{\{1,3,7\}}\delta T)^{-1}T^{\{1,3,7\}}=TT^{\{1,3,7\}}.$$ \vspace*{-0.5\baselineskip}
Then \vspace*{-0.5\baselineskip} $$(\overline{T}B)^*=(TT^{\{1,3,7\}})^*=TT^{\{1,3,7\}}=\overline{T}B,$$ \vspace*{-0.5\baselineskip}
and \vspace*{-1\baselineskip}
\begin{eqnarray*}
\overline{T}B^2&=&TT^{\{1,3,7\}}B\\
&=&TT^{\{1,3,7\}}T^{\{1,3,7\}}(I+\delta TT^{\{1,3,7\}})^{-1}\\
 &=&T^{\{1,3,7\}}(I+\delta TT^{\{1,3,7\}})^{-1}=B.
 \end{eqnarray*}
 It is easy to see $\|B\|\leq \|T^{\{1,3,7\}}\|\cdot\|(I+T^{\{1,3,7\}}\delta T)^{-1}\|$ and  \vspace*{-0.2\baselineskip}$$ \|B-T^{\{1,3,7\}}\|\leq \|T^{\{1,3,7\}}\|\cdot \|(I+T^{\{1,3,7\}}\delta T)^{-1}\|\cdot\|T^{\{1,3,7\}}\delta T\|.$$
 \vspace*{-0.3\baselineskip} This completes the proof.
\par\vspace{0.1in}
If $\|T^{\{1,3,7\}}\delta T\|<1$, then $I+T^{\{1,3,7\}}\delta T: X\rightarrow X$ is invertible. And we can get the following corollary which contains Theorems 2.1, 3.1 and 3.2 in \cite{M}.
 \begin{Corollary}\label{Cor}${\rm\cite{M}}$  Let $T\in C^{n\times n}$ be a matrix
with a  core inverse $T^{\tiny\textcircled{\tiny\#}}\in C^{n\times n}$ and
$\delta T\in C^{n\times n}$ with $\|T^{\tiny\textcircled{\tiny\#}}\delta T\|<1$. If $TT^{\tiny\textcircled{\tiny\#}}\delta T=\delta T$ or $T^{\tiny\textcircled{\tiny\#}} T\delta T= \delta T$, then \vspace*{-0.3\baselineskip}
$$B=T^{\tiny\textcircled{\tiny\#}}(I+\delta TT^{\tiny\textcircled{\tiny\#}})^{-1}=(I+T^{\tiny\textcircled{\tiny\#}}\delta T)^{-1}T^{\tiny\textcircled{\tiny\#}}$$  \vspace*{-0.5\baselineskip} is a  core inverse of $\overline{T}=T+\delta T$, $\overline{T}B=TT^{\tiny\textcircled{\tiny\#}}$ and $$\frac{\|T^{\tiny\textcircled{\tiny\#}}\|}{1+\|T^{\tiny\textcircled{\tiny\#}}\delta T\|}\leq \|B\|\leq\frac{\|T^{\tiny\textcircled{\tiny\#}}\|}{1-\|T^{\tiny\textcircled{\tiny\#}}\delta T\|}.$$
\end{Corollary}
\proof If $T^{\tiny\textcircled{\tiny\#}} T\delta T= \delta T$, then $TT^{\tiny\textcircled{\tiny\#}}\delta T=TT^{\tiny\textcircled{\tiny\#}} T^{\tiny\textcircled{\tiny\#}} T\delta T=T^{\tiny\textcircled{\tiny\#}} T\delta T=\delta T$. It follows from Theorem 2.1 that $B$ is a  core inverse of $\overline{T}$.  This completes the proof.
\begin{Remark} In fact, we can prove that $TT^{\tiny\textcircled{\tiny\#}}\delta T=\delta T$ and $T^{\tiny\textcircled{\tiny\#}} T\delta T= \delta T$ are equivalent.
\end{Remark}
\quad\ Now we can consider the case of the  core inverse in Hilbert space, which is  precisely the $\{1,2,3,6,7\}$-inverse \cite {RDD}.
\begin{Theorem}\label{The2.5} Let $X$ be a Hilbert
 space.   Let $T\in B(X)$
with a  core inverse $T^{\tiny\textcircled{\tiny\#}}\in B(X)$ and
$I+T^{\tiny\textcircled{\tiny\#}}\delta T: X\rightarrow X$ be invertible  with $\delta T\in B(X)$.
Then \vspace*{-0.3\baselineskip}
$$B=T^{\tiny\textcircled{\tiny\#}}(I+\delta TT^{\tiny\textcircled{\tiny\#}})^{-1}=(I+T^{\tiny\textcircled{\tiny\#}}\delta T)^{-1}T^{\tiny\textcircled{\tiny\#}}$$ \vspace*{-0.5\baselineskip}  is a  core inverse of $\overline{T}=T+\delta T$ if and only if
$$R(\overline{T})\subseteq R(T)$$ \vspace*{-0.5\baselineskip}
if and only if \vspace*{-0.5\baselineskip}
$$R(\overline{T})= R(T)$$
 \vspace*{-0.5\baselineskip}if and only if
$$ TT^{\tiny\textcircled{\tiny\#}}\overline{T}= \overline{T} \quad or \quad TT^{\tiny\textcircled{\tiny\#}}\delta T= \delta T. $$
 \vspace*{-0.3\baselineskip}\quad\ In this case, $\overline{T}B=TT^{\tiny\textcircled{\tiny\#}}$, $\|B\|\leq \|T^{\tiny\textcircled{\tiny\#}}\|\cdot \|(I+T^{\tiny\textcircled{\tiny\#}}\delta T)^{-1}\|$ and $$ \|B-T^{\tiny\textcircled{\tiny\#}}\|\leq \|T^{\tiny\textcircled{\tiny\#}}\|\cdot \|(I+T^{\tiny\textcircled{\tiny\#}}\delta T)^{-1}\|\cdot\|T^{\tiny\textcircled{\tiny\#}}\delta T\|.$$
\end{Theorem}
\proof It follows from Theorem 2.1 that we only need to show the ``if'' part. If $R(\overline{T})\subseteq R(T)$, similar to the proof in Theorem 2.1,  we can have that $B$ is a $\{1,3,7\}$-inverse of $\overline{T}$ and \vspace*{-0.5\baselineskip} $$R(\overline{T})= R(T),\quad
\overline{T}=TT^{\tiny\textcircled{\tiny\#}}\overline{T},\quad \overline{T}=T(I+T^{\tiny\textcircled{\tiny\#}}\delta T)\quad {\rm and }\quad \overline{T}B=TT^{\tiny\textcircled{\tiny\#}}.$$ \vspace*{-0.6\baselineskip} Hence \vspace*{-0.8\baselineskip}
$$B\overline{T}B=BTT^{\tiny\textcircled{\tiny\#}}=(I+T^{\tiny\textcircled{\tiny\#}}\delta T)^{-1}T^{\tiny\textcircled{\tiny\#}} TT^{\tiny\textcircled{\tiny\#}}=(I+T^{\tiny\textcircled{\tiny\#}}\delta T)^{-1}T^{\tiny\textcircled{\tiny\#}}=B,$$ \vspace*{-1\baselineskip}
and \vspace*{-0.8\baselineskip}
\begin{eqnarray*}
B\overline{T}^2&=&(I+T^{\tiny\textcircled{\tiny\#}}\delta T)^{-1}T^{\tiny\textcircled{\tiny\#}}(T+\delta T)\overline{T}\\
 &=&(I+T^{\tiny\textcircled{\tiny\#}}\delta T)^{-1}(T^{\tiny\textcircled{\tiny\#}} T-I+I+T^{\tiny\textcircled{\tiny\#}}\delta T)\overline{T}\\
 &=&(I+T^{\tiny\textcircled{\tiny\#}}\delta T)^{-1}(T^{\tiny\textcircled{\tiny\#}} T-I)\overline{T}+\overline{T}\\
 &=&(I+T^{\tiny\textcircled{\tiny\#}}\delta T)^{-1}(T^{\tiny\textcircled{\tiny\#}} T-I)TT^{\tiny\textcircled{\tiny\#}}\overline{T}+\overline{T}\\
 &=&(I+T^{\tiny\textcircled{\tiny\#}}\delta T)^{-1}(T^{\tiny\textcircled{\tiny\#}} TTT^{\tiny\textcircled{\tiny\#}}-TT^{\tiny\textcircled{\tiny\#}})\overline{T}+\overline{T}\\
  &=&(I+T^{\tiny\textcircled{\tiny\#}}\delta T)^{-1}(TT^{\tiny\textcircled{\tiny\#}}-TT^{\tiny\textcircled{\tiny\#}})\overline{T}+\overline{T}=\overline{T}.
 \end{eqnarray*}
Thus $B$ is a $\{2,6\}$-inverse of $\overline{T}$ and so $B$ is the core inverse of $\overline{T}$. This completes the proof.
\par\vspace{0.05in}
 Next, we shall present two examples of $4\times 4$ matrices to illustrate our main results.
  \vspace*{-0.5\baselineskip}
 \begin{Example}\label{Exa3.1}
{\rm Let {\setlength\abovedisplayskip{1pt}
\setlength\belowdisplayskip{1pt}
\[ T = \left[ \begin{array}{rrrr}
1 & 0 & 2 & 4  \\
2 & 1 & -1 & 0  \\
2 & 2 & 0 & 1 \\
1 & -2 & 0 & 2
             \end{array} \right]\; \; \mbox{and} \; \; \; \delta T = \left[ \begin{array}{rrrr}
0 & -1 & 0 & -4 \\
 -2 &  -2 & -2 & 2 \\
-4 & -2 & -4 & 0 \\
4 &  -1 & 4 & 0
       \end{array} \right], \]}
  then \vspace*{-0.5\baselineskip}  {\setlength\abovedisplayskip{1pt}
\setlength\belowdisplayskip{1pt}
\[T^{\tiny\textcircled{\tiny\#}}= \frac{1}{120} \left[ \begin{array}{rrrr}
-30 & 60 & 40 & 10  \\
21 & -18 & 8 & -31  \\
-15 & 30 & 40 & -35 \\
42 & -36 & -24 & 18
             \end{array} \right],\; \; \; \overline{T} = \left[ \begin{array}{rrrr}
1 & -1 & 2 & 0  \\
0 & -1 & -3 & 2  \\
-2 & 0 & -4 & 1 \\
5& -3 & 4 & 2
             \end{array} \right]\]}
and \vspace*{-0.5\baselineskip} {\setlength\abovedisplayskip{1pt}
\setlength\belowdisplayskip{1pt}
\[ R(\overline{T})=R(T)=\left\{   \lambda_1\left[ \begin{array}{rrrr}
1  \\
2  \\
2 \\
1
             \end{array} \right]+\lambda_2\left[ \begin{array}{rrrr}
0  \\
1  \\
2 \\
-2
             \end{array} \right]+\lambda_3\left[ \begin{array}{rrrr}
2  \\
-1  \\
0 \\
0
             \end{array} \right]:\ \lambda_1,\lambda_2,\lambda_3\in C \right\}.\]}
 By Theorem 2.1,  $(I+T^{\tiny\textcircled{\tiny\#}}\delta T)^{-1}T^{\tiny\textcircled{\tiny\#}}$ is the core inverse of $\overline{T}$. It should be pointed out that  a direct computation can also give
 {\setlength\abovedisplayskip{1pt}
\setlength\belowdisplayskip{1pt}
\[ \overline{T}^{\tiny\textcircled{\tiny\#}}=(I+T^{\tiny\textcircled{\tiny\#}}\delta T)^{-1}T^{\tiny\textcircled{\tiny\#}}=\frac{1}{90} \left[ \begin{array}{rrrr}
30 & 0 & 10 & 10  \\
-201 & 18 & -88 & 11  \\
-75 & 0 & -40 & 5 \\
-222 & 36 & -86 & 22
             \end{array} \right].\]}}
\end{Example}
 \begin{Example}\label{Exa3.2}\
{\rm Let $T$ be the same matrix in Example 2.1 and {\setlength\abovedisplayskip{1pt}
\par
\setlength\belowdisplayskip{1pt}
 \[ \delta T = \left[ \begin{array}{rrrr}
1 & 0 & -2 & -2 \\
 -2 &  1 & 1 & -1 \\
-2 & -2 & 2 & 2 \\
-1 &  2 & 0 & -2
       \end{array} \right], \]}
  then   {\setlength\abovedisplayskip{1pt}
\setlength\belowdisplayskip{1pt}
\[\overline{T} = \left[ \begin{array}{rrrr}
2 & 0 & 0 & 2  \\
0 & 2 & 0 & -1  \\
0 & 0 & 2 & 3 \\
0 & 0 & 0 & 0
             \end{array} \right],\; \;\; \; \left[ \begin{array}{rrrr}
1  \\
2  \\
2 \\
1 \\
             \end{array} \right]\in R(T)\; \; \; \mbox{and} \; \; \;   \left[ \begin{array}{rrrr}
1  \\
2  \\
2 \\
1 \\
\end{array} \right]\notin R(\overline{T}),\]}
\noindent which implies $R(\overline{T})\neq R(T)$. Hence by Theorem 2.1, $(I+T^{\tiny\textcircled{\tiny\#}}\delta T)^{-1}T^{\tiny\textcircled{\tiny\#}}$ is not the core inverse of $\overline{T}$. In addition, we can also verify $\overline{T}^{\tiny\textcircled{\tiny\#}}\neq(I+T^{\tiny\textcircled{\tiny\#}}\delta T)^{-1}T^{\tiny\textcircled{\tiny\#}}$ directly by
\par
 {\setlength\abovedisplayskip{1pt}
\setlength\belowdisplayskip{1pt}
\[ \overline{T}^{\tiny\textcircled{\tiny\#}}=\frac{1}{2}\left[ \begin{array}{rrrr}
1 & 0 & 0 & 0  \\
0 & 1 & 0 & 0  \\
0 & 0 & 1 & 0 \\
0 & 0 & 0 & 0
             \end{array} \right].\]} and\vspace*{-0.5\baselineskip}
 {\setlength\abovedisplayskip{1pt}
\setlength\belowdisplayskip{1pt}
\[ (I+T^{\tiny\textcircled{\tiny\#}}\delta T)^{-1}T^{\tiny\textcircled{\tiny\#}}=\frac{1}{8} \left[ \begin{array}{rrrr}
6 & 4 & -4 & 22  \\
-1 & 2 & 2 & -1  \\
3 & 6 & -2 & 19 \\
-2 & -4 & 4 & -18
             \end{array} \right].\]}}
\end{Example}\vspace*{-0.5\baselineskip}
\begin{Remark} According to $N(\overline{T})=N(T)$ and ${\rm Rank}\ \overline{T}={\rm Rank}\ T$ in Example 2.2,  we can conclude that the null space-preserving and rank-preserving perturbation can not guarantee that the
core inverse of the perturbed operator has the simplest possible expression.
\end{Remark}\vspace*{-1\baselineskip}
\section{Conclusions}\vspace*{-0.6\baselineskip}
\quad\ As we have seen, the core inverse $T^{\tiny\textcircled{\tiny\#}}$ is closely related to the orthogonal projector on the range $R(T)$. This is the principal reason that we can find the characterization for the
core inverse of the perturbed operator to have the simplest possible expression. How about the
expressions for the core inverse under the null space-preserving perturbation, rank-preserving perturbation, more generally, the stable perturbation or acute perturbation \cite{HZJ,We,WeD,WeW}? We hope to
solve these more and general cases in the future.\vspace*{-0.6\baselineskip}

\end{document}